\documentclass[11pt, twoside]{article}
\usepackage{latexsym}
\usepackage{amsmath}
\usepackage{amssymb}
\usepackage[all]{xy}
\usepackage{amsfonts}
\usepackage{verbatim}
\usepackage{amsthm}
\usepackage{mathrsfs}
\usepackage{epsfig}
\usepackage{xy}
\usepackage{array}
\usepackage{stmaryrd}
\usepackage{graphicx,color}
\usepackage{xcolor}
\usepackage{tikz}
\usetikzlibrary{arrows,calc}
\usepackage{etex}
\usepackage{mathdots}
\usepackage{float}
\usepackage{graphics}
\usepackage{pdflscape}

\usepackage{anysize,hyperref}
\input xypic
\xyoption{all}

\usepackage[perpage,symbol]{footmisc}
\usepackage{setspace}
\def\dr{\ar@{->}[r]}

\begin{document}
\baselineskip=15pt
\title{\Large{\bf Square commutative groups}}
\medskip
\author{Weicai Wu\footnote{Weicai Wu was supported by the Project of Improving the Basic Scientific Research Ability of Young and Middle-aged College Teachers in Guangxi (Grant No: 2024KY0069).}, Mingxuan Yang, Yangbo Zhou and Chao Rong
}

\date{}

\maketitle
\def\blue{\color{blue}}
\def\red{\color{red}}

\newtheorem{Theorem}{Theorem}[section]
\newtheorem{Lemma}[Theorem]{Lemma}
\newtheorem{Corollary}[Theorem]{Corollary}
\newtheorem{Proposition}[Theorem]{Proposition}
\newtheorem{Conjecture}{Conjecture}
\theoremstyle{Definition}
\newtheorem{Definition}[Theorem]{Definition}
\newtheorem{Question}[Theorem]{Question}
\newtheorem{Remark}[Theorem]{Remark}
\newtheorem{Remark*}[]{Remark}
\newtheorem{Example}[Theorem]{Example}
\newtheorem{Example*}[]{Example}
\newtheorem{Condition}[Theorem]{Condition}
\newtheorem{Condition*}[]{Condition}
\newtheorem{Construction}[Theorem]{Construction}
\newtheorem{Construction*}[]{Construction}
\newtheorem{Assumption}[Theorem]{Assumption}
\newtheorem{Assumption*}[]{Assumption}

\baselineskip=17pt
\parindent=0.5cm
\vspace{-6mm}

\begin{abstract}
\baselineskip=16pt
In this paper we first give a necessary and sufficient  condition for a group $G$ generated by $n$ elements to be a square commutative group and prove $G$ is a square commutative group if and only if $\widehat{G}$ is an abelian group, then we give conditions for a group generated by two elements, with additional conditions, to be a square commutative group.\\[0.2cm]
\textbf{Key words:} Square commutative groups;  Generators; Baumslag-Solitar groups.\\[0.1cm]
\textbf{ 2020 Mathematics Subject Classification:} 20E05; 20F05.
\end{abstract}

\pagestyle{myheadings}
\markboth{\rightline {\scriptsize   Weicai Wu, Mingxuan Yang, Yangbo Zhou and Chao Rong}}
         {\leftline{\scriptsize Square commutative groups}}

\section{Introduction}
The structure of Nichols algebras over finite abelian groups is well-known in \cite{H}. It has been shown in \cite {HS} that if a Nichols algebra over non-abelian groups is finite dimensional, then concrete elements of the conjugacy classes must satisfy the square commutative law. In \cite {AFGV}, an equivalence between the rack operation and square commutative has been found, and a large type of infinite-dimensional Nichols algebras was eliminated by using rack theory. It is a quite powerful sufficient condition for infinite-dimensionality, at least if the rack is sufficiently large and far from being abelian. But unfortunately, rack theory fails in the classification of Nichols algebras over square commutative groups.
Since all braided vector spaces of diagonal type can be obtained from Yetter-Drinfeld (YD) modules over abelian groups, we can consider YD modules on the square commutative group, and the first step requires giving a decomposition of the square commutative group.
On the other hand, the square commutative law is an important arithmetic rule in algebra, and similar to the commutative law, we can also define the square commutative group and square commutative algebra. One defines the square commutative group in \cite {W}, and proves that the groups of order less than $12$, only the dihedral group $D_{6}$ and $D_{10}$ are not square commutative groups. In this paper, we obtain a necessary and sufficient condition for a group generated by $n$ elements to be a square commutative group and establish the relationship between square commutative groups and abelian groups. Our results bring some new light to the classification of Nichols algebras over finite square commutative groups.

Any group is isomorphic to the quotient group of a free group, and any finitely generated group is isomorphic to the quotient group of a finitely generated free group. Using the free group we can give a presentation of the group, e.g. the dihedral group $D_{2n}$ can be expressed in the form $\{1,a,\ldots,a^{n-1},b,ab,\ldots,a^{n-1}b\}$, where $a^{n}=b^{2}=1,aba=b$. We can view it as the free groups generated by $a,b$ with the additional constraints, here we call $a,b$ the generators of $D_{2n}$, the constraints $a^{n}=b^{2}=1,aba=b$ are called the generating relations of $a,b$, and $a^{n}=b^{2}=1$ and $aba=b$ are also called the two relational formulas satisfied by $a,b$, $D _{2n}$ can also be expressed as
$\langle a,b\mid a^{n}=b^{2}=1,aba=b\rangle$.
In this paper, we consider the square commutative groups $G$ generated by $a,b$ and $T_{n}=\{a_{1},a_{2},\ldots,a_{n}\mid n\geq3\}$, respectively. Set $X\cdot Y:=\{xy \mid   x\in X,y\in Y\}$ and $X^{[2]}:=\{x^{2}\mid x\in X\}$ if $X,Y\subset G$. $Z_{G}$ denotes the center of group $G$. Let $Z_{G}^{2}:=\{d\in Z_{G} \mid   d^{2}=e\}$ and $\widehat{G}:=G/Z_{G}^{2}$. These are our four main results.

{\bf Theorem 2.6} {\it $G$ is a square commutative group if and only if $a,b$ satisfy the relation $b^{2}a=ab^{2}, a^{2}b=ba^{2}, (ab)^{2}=(ba)^{2}$, where $G\subset C_{2}\cdot Z_{G}$.}

{\bf Theorem 3.4} {\it Assume that $n\geq3$, then $G$ is a square commutative group if and only if $x_{1}x_{2}^{2}=x_{2}^{2}x_{1}$, $x_{1}(x_{2}x_{3})^{2}=(x_{2}x_{3})^{2}x_{1}$,
$(x_{1}x_{2})^{2}=(x_{2}x_{1})^{2}$, where any $x_{1},x_{2},x_{3}$ are mutually exclusive elements in $T_{n}$, in which case we have $G\subset C_{n}\cdot Z_{G}$.}

{\bf Theorem 4.4} {\it $G$ is a square commutative group if and only if $\widehat{G}$ is an abelian group. Meanwhile, we have $xy\sim yx$ for any $x,y\in G$.}

{\bf Theorem 4.8} {\it Assume that $G$ is a group. Then $G^{[2]}\subset Z_{G}$ if and only if $G$ is a square commutative group.}

This paper is organized as follows. In section 1 we recall some basic results on groups and fix the notation.
In section 2 and section 3 we present the necessary and sufficient condition for $G$ to be the square commutative groups generated by $n$ elements and give the center of the square commutative group. In section 4 we prove $G$ is a square commutative group if and only if $\widehat{G}$ is an abelian group. In section 5
we give the conditions for $G$ to be a square commutative group generated by $a,b$ with defining relations $a^{p}b=ba^{q}$ and some
additional relations.

Throughout, $\mathbb Z:=\{x \mid   x \hbox { is an integer}\}$, $\mathbb N :=\{x \mid   x \hbox { is an integer},  x>0\}$, $\mathbb N_{0} :=\{x \mid   x \hbox { is an integer},  x\geq0\}$, $\mathbb Z^{*} :=\{x \mid   x \hbox { is an integer},  x\neq0\}$, $\mathbb R := \{ x \mid x \hbox { is a real number}\}$. $\mid x\mid $ denotes the order of group $x$ or group element $x$.

\section{Conditions for square commutative groups generated by two elements}
Without any specific statement, the following $G$ in this section is the group generated by $a,b$ with unit element $e$.

\begin {Definition}\label {1.0}(See \cite [Def 1]{W})
A group $G$ is called the square commutative group if $(xy)^{2}=(yx)^{2}$ for all $x,y\in G$, otherwise, $G$ is called the non square commutative group.
\end {Definition}

\begin {Proposition} \label {1.1} {\rm (i)} Assuming that $G$ is a square commutative group. For all $x,y\in G$, then $(xyx)^{m}y^{n}=y^{n}(xyx)^{m}$ holds for any
$m,n\in \mathbb Z$.

{\rm (ii)} Assuming that $G$ is a free group generated by $S=\{a_{1},a_{2},\cdots, a_{n}\}$, then $G$ is a abelian group if and only if $a_{i}a_{j}=a_{j}a_{i}$ holds for any $1\leq i,j\leq n$.
\end {Proposition}

\noindent {\it Proof.} {\rm (i)} can be proved by induction. {\rm (ii)} is obvious.
\hfill $\Box$

If $G$ is a group generated by $a,b$, then $G$ is a abelian group if and only if $ab=ba$ by Proposition \ref {1.1} {\rm (ii)}. But the following conclusion no longer holds: $G$ is a square commutative group if and only if $(ab)^{2}=(ba)^{2}$.

\begin {Example} \label {1.1'} $G=\langle a,b\mid a^{3}=b^{2}=e,ab=ba^{2}\rangle=\{e,a,a^{2},b,ba,ba^{2}\}$, it is clear $(ab)^{2}=(ba)^{2}$, but $G$ is a non square commutative group since $(bba)^{2}\neq(bab)^{2}$.
\end {Example}

\begin {Lemma} \label {1.2} If $G$ is a square commutative group, then $(ab)^{2}=(ba)^{2}, a^{2}b=ba^{2}, b^{2}a=ab^{2}$.
\end {Lemma}

\noindent {\it Proof.} Since $G$ is a square commutative group, then $abab=baba,(bab)^{2}=(abb)^{2}$, i.e., $bab=ababa^{-1}, bab^{2}ab=ab^{2}ab^{2}$, and $ababa^{-1}ba=(bab)ba=ab(bab)=abababa^{-1}$, simplifying, $ba^{2}=a^{2}b$. Similarly, $ab^{2}=b^{2}a$.
\hfill $\Box$

Set $C_{2}:=\{e,a,b,ab\}$.

\begin {Proposition} \label {1.3'} If $a,b$ satisfy the relation $ab^{2}=b^{2}a,ba^{2}=a^{2}b,(ab)^{2}=(ba)^{2}$, then
$a^{p},(ab)^{p}$, $b^{p}\in Z_{G}$, where $\frac{p}{2}\in \mathbb Z$, and $x^{2}\in Z_{G}$ for any $x\in C_{2}$.
\end {Proposition}
\noindent {\it Proof.} Since $ab^{2}=b^{2}a,ba^{2}=a^{2}b,(ab)^{2}=(ba)^{2}$,
then $xy^{2}=y^{2}x,yx^{2}=x^{2}y,(xy)^{2}=(yx)^{2}$, for any $x\in\{a,a^{-1}\},
y\in\{b,b^{-1}\}$. The remaining results are obvious.
\hfill $\Box$

\begin {Proposition} \label {1.3} If $a,b$ satisfy the relation $ab^{2}=b^{2}a,ba^{2}=a^{2}b,(ab)^{2}=(ba)^{2}$, then $G\subset C_{2}\cdot Z_{G}$ and $G$ is a square commutative group.
\end {Proposition}

\noindent {\it Proof.} Since $e\in C_{2}\cdot Z_{G}$, it is only necessary to show that $u_{1}^{n_{1}}u_{2}^{n_{2}}\cdots u_{k}^{n_{k}}\in C_{2}\cdot Z_{G}$ for any $1\leq j\leq k, u_{j}\in \{a,b\},u_{l}\neq u_{l+1},1\leq l\leq k-1,n_{j}\in\mathbb Z^{*}$ are valid. Use induction on $k$. When $k=1$, $u_{1}^{n_{1}}\in C_{2}\cdot Z_{G}$ is obvious.
When $k>1$, assuming that the conclusion holds for all $s\leq k-1$, consider $u_{1}^{n_{1}}u_{2}^{n_{2}}\cdots u_{k}^{n_{k}}$, if
$n_{k}$ is even, then $u_{k}^{n_{k}}\in Z_{G}$, $u_{1}^{n_{1}}u_{2}^{n_{2}}\cdots u_{k}^{n_{k}}\in C_{2}\cdot Z_{G}\cdot Z_{G}=C_{2}\cdot Z_{G}$. If $n_{k}$ is odd, then $C_{2}\cdot \{a^{n_{k}}\}=$
$\{a(a^{n_{k}-1}),e(a^{n_{k}+1}),
ab(a^{n_{k}+1}b^{2}(ab)^{-2}),b(a^{n_{k}+3}b^{2}(ab)^{-2})\}\subset C_{2}\cdot Z_{G}$,
$C_{2}\cdot \{b^{n_{k}}\}=\{b(b^{n_{k}-1}),$
$ab(b^{n_{k}-1}),e(b^{n_{k}+1}),
a(b^{n_{k}+1})\}\subset C_{2}\cdot Z_{G}$.
Summarizing, we have $u_{1}^{n_{1}}u_{2}^{n_{2}}\cdots u_{k}^{n_{k}}\in
C_{2}\cdot Z_{G}\cdot \{u_{k}^{n_{k}}\}=Z_{G}\cdot C_{2}\cdot \{u_{k}^{n_{k}}\}\subset Z_{G}\cdot C_{2}\cdot Z_{G}=C_{2}\cdot Z_{G}\cdot Z_{G}=C_{2}\cdot Z_{G}$. Since all elements in $Z_{G}$ are commutative with any element of $G$, a simple check shows that any two elements in $C_{2}$ are square commutative. For all $g_{1},g_{2}\in G$, $g_{1}=c_{1}d_{1},g_{2}=c_{2}d_{2},c_{1},c_{2}\in C_{2},d_{1},d_{2}\in Z_{G}$,  so $(g_{1}g_{2})^{2}=(c_{1}c_{2})^{2}d_{1}^{2}d_{2}^{2}
=(c_{2}c_{1})^{2}d_{1}^{2}d_{2}^{2}=(g_{2}g_{1})^{2}$.
\hfill $\Box$

\begin {Theorem} \label {1.4}  $G$ is a square commutative group if and only if $a,b$ satisfy the relation $(ab)^{2}=(ba)^{2}, a^{2}b=ba^{2}, b^{2}a=ab^{2}$, where $G\subset C_{2}\cdot Z_{G}$.
\end {Theorem}
\noindent {\it Proof.} The conclusion follows from Lemma \ref {1.2} and Proposition \ref {1.3}.
\hfill $\Box$

\begin {Corollary} \label {1.5'}
If $G$ is a square commutative group, then $x=a^{h_{a}}b^{h_{b}}(ab)^{2\lambda_{a,b}}$ for any $x\in G$ with $h_{a},h_{b}\in\mathbb Z,\lambda_{a,b}\in \{0,1\}$.
\end {Corollary}

\begin {Corollary} \label {1.5}
Assume that $b^{2}=e$, then $G$ is a square commutative group if and only if $a,b$ satisfy the relations $(ab)^{2}=(ba)^{2}, a^{2}b=ba^{2}$.
\end {Corollary}

\begin {Theorem} \label {1.6} Suppose $a,b$ satisfy the relation $a^{n}=e$ or $b^{n}=e$, $n$ is odd, then $G$ is a square commutative group if and only if $G$ is an abelian group.
\end {Theorem}

\noindent {\it Proof.} If $a^{n}=e$, let $n=2l-1,l\in \mathbb N$, since $G$ is a square commutative group, then $a^{2}b=ba^{2}, a^{2l}b=ba^{2l}$, thus $a^{n+1}b=ba^{n+1},ab=ba$, the case $b^{n}=e$ is similar.
\hfill $\Box$

\begin {Corollary} \label {1.7} Suppose $a,b$ satisfy the relation $a^{n}=e$ or $b^{n}=e$, $n$ is odd, then $(ab)^{2}=(ba)^{2}, a^{2}b=ba^{2}, b^{2}a=ab^{2}$ if and only if $ab=ba$.
\end {Corollary}

If $a,b$ satisfy the relation $a^{2}=b^{2}=e$, then $G$ is a square commutative group if and only if $a,b$ satisfy $(ab)^{2}=(ba)^{2}$, then $C_{2}=\{e,a,b,ab\}, Z_{G}=\{e,(ab)^{2}\}$, and
$G=C_{2}\cdot Z_{G}$, i.e., $G$ is isomorphic to the dihedral group $D_{8}$, which is the only square commutative group in the dihedral group (see \cite[Proposition 1]{W}).

\section{Square commutative groups generated by $n$ elements
with $n\geq3$}
Let $T_{n}=\{a_{1},a_{2},\ldots,a_{n}\},n\geq3$, without any specific statement the following $G$ in this section is the group generated by $T_{n}$ with unit element $e$. $\mathbb S_{m}$ denotes symmetric group of $m$ elements.

\begin {Lemma} \label {5.1} Assume that $G$ is a square commutative group, then $x_{1}x_{2}^{2}=x_{2}^{2}x_{1}$, $x_{1}(x_{2}x_{3})^{2}=(x_{2}x_{3})^{2}x_{1}$,
$(x_{1}x_{2})^{2}=(x_{2}x_{1})^{2}$, where any $x_{1},x_{2},x_{3}$ are mutually exclusive elements in $T_{n}$.
\end {Lemma}

\noindent {\it Proof.} Let $x_{2,m}:=x_{2}\cdots x_{m}$. Since $G$ is a square commutative group, then $x_{2,m}x_{1}x_{2,m}x_{1}=x_{1}x_{2,m}x_{1}x_{2,m}$,
$(x_{1}x_{2,m}x_{1})^{2}=(x_{2,m}x_{1}x_{1})^{2}$, i.e., $x_{1}x_{2,m}x_{1}=x_{2,m}x_{1}x_{2,m}x_{1}x_{2,m}^{-1}$,
$x_{1}x_{2,m}x_{1}^{2}x_{2,m}x_{1}$
$=x_{2,m}x_{1}^{2}x_{2,m}x_{1}^{2}$, and
$x_{2,m}x_{1}x_{2,m}x_{1}x_{2,m}^{-1}x_{1}x_{2,m}=
(x_{1}x_{2,m}x_{1})x_{1}x_{2,m}
=x_{2,m}x_{1}(x_{1}x_{2,m}x_{1})=
x_{2,m}x_{1}x_{2,m}x_{1}x_{2,m}x_{1}x_{2,m}^{-1}$, simplifying, $x_{1}(x_{2,m})^{2}=(x_{2,m})^{2}x_{1}$, so $x_{1}x_{2}^{2}=x_{2}^{2}x_{1}$,
$x_{1}(x_{2}x_{3})^{2}=(x_{2}x_{3})^{2}x_{1}$.
\hfill $\Box$

Set $C_{n}:=\{e,a_{i_{1}}a_{i_{2}}\cdots a_{i_{j}}\mid 1\leq i_{1}<\cdots<i_{j}\leq n,1\leq j\leq n\}$. It is clear $G\cdot G=G$.

\begin {Proposition} \label {5.3} Assume that $1\leq i_{1},i_{2},\cdots ,i_{m}\leq n,m\in \mathbb N$. If $x_{1}x_{2}^{2}=x_{2}^{2}x_{1}, x_{1}(x_{2}x_{3})^{2}=(x_{2}x_{3})^{2}x_{1},
(x_{1}x_{2})^{2}=(x_{2}x_{1})^{2}$, where any $x_{1}$, $x_{2}$, $x_{3}$ are mutually exclusive elements in $T_{n}$. Then

{\rm (i)} $x_{2}x_{1}=x_{1}x_{2}x_{1}^{2}x_{2}^{2}(x_{1}x_{2})^{-2}$
and $(x_{1}x_{2})^{4}=x_{1}^{4}x_{2}^{4}$,
$\forall\ x_{1},x_{2}\in T_{n}$.

{\rm (ii)} there exists $d\in Z_{G}^{2}$
such that $a_{i_{\tau(1)}}a_{i_{\tau(2)}}\cdots a_{i_{\tau(m)}}=a_{i_{1}}a_{i_{2}}\cdots a_{i_{m}}d$ for any $\tau\in \mathbb S_{m}$.

{\rm (iii)} any two elements in $C_{n}$ are square commutative.

{\rm (iv)} $a_{i}^{p_{i}},(a_{j}a_{k})^{q_{jk}}\in Z_{G}$, where $\frac{p_{i}}{2},\frac{q_{jk}}{2}\in \mathbb Z,1\leq i\leq n,1\leq j<k\leq n$, and $C_{n}^{[2]},G^{[2]}\subset Z_{G}$.
\end {Proposition}
\noindent {\it Proof.} {\rm (i)} If $x_{1}=x_{2}$, it is clear. If
$x_{1}\neq x_{2}$, then $x_{2}x_{1}
=(x_{1}x_{2})^{-2}(x_{1}x_{2}x_{1}x_{2})x_{2}x_{1}
=x_{1}x_{2}x_{1}^{2}x_{2}^{2}(x_{1}x_{2})^{-2}$ and $(x_{1}x_{2})^{4}=(x_{1}x_{2})^{2}(x_{2}x_{1})^{2}=x_{1}^{4}x_{2}^{4}$.

{\rm (ii)} If $\tau=(1)$, then $d=e$. Now assume $\tau\neq(1)$. We show $a_{i_{\tau(1)}}a_{i_{\tau(2)}}\cdots a_{i_{\tau(m)}}=a_{i_{1}}a_{i_{2}}\cdots a_{i_{m}}d$ by induction on $m$. If $m=2$, $a_{i_{\tau(1)}}a_{i_{\tau(2)}}=a_{i_{2}}a_{i_{1}}=
a_{i_{1}}a_{i_{2}}a_{i_{1}}^{2}a_{i_{2}}^{2}(a_{i_{1}}a_{i_{2}})^{-2}$, $d^{2}=
a_{i_{1}}^{4}a_{i_{2}}^{4}(a_{i_{1}}a_{i_{2}})^{-4}=e$ by {\rm (i)}. We assume $m>2$ from now on. If $\tau(m)=m$, The claim holds by the induction hypothesis and $\tau\in \mathbb S_{m-1}$. If $\tau(k)=m$, $1\leq k<m$, then $
\{i_{\tau(1)},\ldots, i_{\tau(k-1)},i_{\tau(k+1)},\ldots, i_{\tau(m)}\}=\{i_{1},\ldots,i_{m-1}\}$.

$a_{i_{\tau(1)}}\cdots (a_{i_{\tau(k)}}a_{i_{\tau(k+1)}})\cdots a_{i_{\tau(m)}}$

$=a_{i_{\tau(1)}}\cdots (a_{i_{\tau(k+1)}}a_{i_{\tau(k)}}a_{i_{\tau(k+1)}}^{2}a_{i_{\tau(k)}}^{2}
(a_{i_{\tau(k+1)}}a_{i_{\tau(k)}})^{-2})\cdots a_{i_{\tau(m)}}$

$=a_{i_{\tau(1)}}\cdots(a_{i_{\tau(k+1)}}a_{i_{\tau(k)}})\cdots a_{i_{\tau(m)}}a_{i_{\tau(k+1)}}^{2}a_{i_{\tau(k)}}^{2}
(a_{i_{\tau(k+1)}}a_{i_{\tau(k)}})^{-2}$

$=a_{i_{\tau(1)}}\cdots a_{i_{\tau(k-1)}}a_{i_{\tau(k+1)}}\cdots a_{i_{\tau(m)}}a_{i_{\tau(k)}}\prod\limits_{j=k+1}^{m}a_{i_{\tau(j)}}^{2}
(a_{i_{\tau(j)}}a_{i_{\tau(k)}})^{-2}a_{i_{\tau(k)}}^{2(m-k)}
$

$=a_{i_{\tau(1)}}\cdots a_{i_{\tau(k-1)}}a_{i_{\tau(k+1)}}\cdots a_{i_{\tau(m)}}a_{i_{m}}\prod\limits_{j=k+1}^{m}a_{i_{\tau(j)}}^{2}
(a_{i_{\tau(j)}}a_{i_{m}})^{-2}a_{i_{m}}^{2(m-k)}
$

$=a_{i_{\tau_{1}(1)}}\cdots a_{i_{\tau_{1}(k-1)}}a_{i_{\tau_{1}(k)}}\cdots a_{i_{\tau_{1}(m-1)}}a_{i_{m}}\prod\limits_{j=k+1}^{m}a_{i_{\tau(j)}}^{2}
(a_{i_{\tau(j)}}a_{i_{m}})^{-2}a_{i_{m}}^{2(m-k)}
$ ( $\tau_{1}\in \mathbb S_{m-1}$ )

$=a_{i_{1}}a_{i_{2}}\cdots a_{i_{m-1}}d_{1}a_{i_{m}}\prod\limits_{j=k+1}^{m}a_{i_{\tau(j)}}^{2}
(a_{i_{\tau(j)}}a_{i_{m}})^{-2}a_{i_{m}}^{2(m-k)}
$( by inductive hypothesis )

$d^{2}=d_{1}^{2}\prod\limits_{j=k+1}^{m}a_{i_{\tau(j)}}^{4}
(a_{i_{\tau(j)}}a_{i_{\tau(k)}})^{-4}a_{i_{\tau(k)}}^{4(m-k)}$( by $d_{1}\in Z_{G}^{2}$ )

$=\prod\limits_{j=k+1}^{m}a_{i_{\tau(j)}}^{4}
(a_{i_{\tau(j)}}a_{i_{\tau(k)}})^{-2}
(a_{i_{\tau(k)}}a_{i_{\tau(j)}})^{-2}a_{i_{\tau(k)}}^{4(m-k)}$

$=\prod\limits_{j=k+1}^{m}(a_{i_{\tau(j)}}^{4}
a_{i_{\tau(j)}}^{-4}a_{i_{\tau(k)}}^{-4})
a_{i_{\tau(k)}}^{4(m-k)}=e$.

{\rm (iii)} For $\forall\ a_{i_{1}}\cdots a_{i_{j}},a_{l_{1}}\cdots a_{l_{k}}\in C_{n}$, $ (a_{i_{1}}\cdots a_{i_{j}})(a_{l_{1}}\cdots a_{l_{k}})=(a_{l_{1}}\cdots a_{l_{k}})(a_{i_{1}}\cdots a_{i_{j}})d$ and $d^{2}=e$ by {\rm (ii)}, then $(a_{i_{1}}\cdots a_{i_{j}}a_{l_{1}}\cdots a_{l_{k}})^{2}=(a_{l_{1}}\cdots a_{l_{k}}a_{i_{1}}\cdots a_{i_{j}})^{2}$.

{\rm (iv)} It is clear.
\hfill $\Box$

\begin {Proposition} \label {5.4} If $x_{1}x_{2}^{2}=x_{2}^{2}x_{1}$, $x_{1}(x_{2}x_{3})^{2}=(x_{2}x_{3})^{2}x_{1}$,
$(x_{1}x_{2})^{2}=(x_{2}x_{1})^{2}$, where any $x_{1},x_{2},x_{3}$ are mutually exclusive elements in $T_{n}$, then

{\rm (i)} $C_{n}\cdot T_{n}\subset C_{n}\cdot Z_{G}$

{\rm (ii)} $G\subset C_{n}\cdot Z_{G}$, and $G$ is a square commutative group.
\end {Proposition}

\noindent {\it Proof.} {\rm (i)} Consider $a_{i_{1}}a_{i_{2}}\cdots a_{i_{j}}a_{t}$. If $t\notin \{i_{1},i_{2},\ldots,i_{j}\}$, then there exist $l_{1},l_{2},\ldots,l_{j+1}$ such that $a_{i_{1}}a_{i_{2}}\cdots a_{i_{j}}a_{t}=a_{l_{1}}a_{l_{2}}\cdots a_{l_{j+1}}d\in C_{n}\cdot Z_{G}$  with $\{l_{1},l_{2},\ldots,l_{j+1}\}=\{i_{1},i_{2},\ldots,i_{j},t\}$ and $l_{1}<l_{2}<\ldots<l_{j+1},d^{2}=e$ by Proposition \ref {5.3}{\rm (ii)}. If $t\in \{i_{1},i_{2},\ldots,i_{j}\}$, then there exist $i_{p}=t$ with $1\leq p\leq j$, then $a_{i_{1}}a_{i_{2}}\cdots a_{i_{j}}a_{t}=a_{i_{1}}a_{i_{2}}\cdots a_{i_{p-1}}a_{i_{p+1}}\cdots a_{i_{j}}a_{t}^{2}d_{1}\in C_{n}\cdot Z_{G}$ and $d_{1}^{2}=e$ by Proposition \ref {5.3}{\rm (ii)}.

{\rm (ii)} Since $e\in C_{n}\cdot Z_{G}$, it is only necessary to show that $u_{1}^{n_{1}}u_{2}^{n_{2}}\cdots u_{k}^{n_{k}}\in C_{n}\cdot Z_{G}$ for any $1\leq j\leq k, u_{j}\in T_{n},u_{l}\neq u_{l+1},1\leq l\leq k-1,n_{j}\in\mathbb Z^{*}$ are valid. Use induction on $k$. When $k=1$, $u_{1}^{n_{1}}\in C_{n}\cdot Z_{G}$ is obvious. When $k>1$, assuming that the conclusion holds for all $s\leq k-1$, consider $u_{1}^{n_{1}}u_{2}^{n_{2}}\cdots u_{k}^{n_{k}}$, if
$n_{k}$ is even, then $u_{k}^{n_{k}}\in Z_{G}$, $u_{1}^{n_{1}}u_{1}^{n_{2}}\cdots u_{k}^{n_{k}}\in C_{n}\cdot Z_{G}\cdot Z_{G}=C_{n}\cdot Z_{G}$. If $n_{k}$ is an odd, then $u_{1}^{n_{1}}u_{2}^{n_{2}}\cdots u_{k}u_{k}^{n_{k}-1}\in C_{n}\cdot Z_{G}\cdot T_{n}\cdot Z_{G}=C_{n}\cdot T_{n}\cdot Z_{G}\cdot Z_{G}\subset C_{n}\cdot Z_{G}$ by {\rm (i)}.
Summarizing, we have $u_{1}^{n_{1}}u_{2}^{n_{2}}\cdots u_{k}^{n_{k}}\in
C_{n}\cdot Z_{G}\cdot \{u_{k}^{n_{k}}\}=Z_{G}\cdot C_{n}\cdot \{u_{k}^{n_{k}}\}\subset Z_{G}\cdot C_{n}\cdot Z_{G}=C_{n}\cdot Z_{G}\cdot Z_{G}=C_{n}\cdot Z_{G}$. Since all elements in $Z_{G}$ are commutative with any element of $G$ and any two elements in $C_{n}$ are square commutative by Proposition \ref {5.3}{\rm (iii)}. For all $g_{1},g_{2}\in G$, $g_{1}=c_{1}d_{1},g_{2}=c_{2}d_{2},c_{1},c_{2}\in C_{n},d_{1},d_{2}\in Z_{G}$,  so $(g_{1}g_{2})^{2}=(c_{1}c_{2})^{2}d_{1}^{2}d_{2}^{2}=(g_{2}g_{1})^{2}$.
\hfill $\Box$

\begin {Theorem} \label {5.5} Assume that $n\geq3$, then $G$ is a square commutative group if and only if $x_{1}x_{2}^{2}=x_{2}^{2}x_{1}$, $x_{1}(x_{2}x_{3})^{2}=(x_{2}x_{3})^{2}x_{1}$,
$(x_{1}x_{2})^{2}=(x_{2}x_{1})^{2}$, where any $x_{1},x_{2},x_{3}$ are mutually exclusive elements in $T_{n}$, in which case we have $G\subset C_{n}\cdot Z_{G}$.
\end {Theorem}
\noindent {\it Proof.} The conclusion follows from Lemma \ref {5.1} and Proposition \ref {5.4}.
\hfill $\Box$

\begin {Corollary} \label {5.5'}
If $G$ is a square commutative group, then $x=a_{1}^{h_{1}}a_{2}^{h_{2}}\cdots a_{n}^{h_{n}}\prod\limits_{1\leq i<j\leq n}(a_{i}a_{j})^{2\lambda_{i,j}}$ for any $x\in G$ with $h_{k}\in\mathbb Z,1\leq k\leq n,\lambda_{i,j}\in \{0,1\}$.
\end {Corollary}

\section{Relationship between square commutative groups and abelian groups}
In this section we prove that $G$ is a square commutative group if and only if $\widehat{G}$ is an abelian group.

\begin {Proposition} \label {5.6} Assume that $G$ is a square commutative group, then $G/Z_{G}$ is an abelian group.
\end {Proposition}
\noindent {\it Proof.} For any $x,y\in G/Z_{G}$, we have $x=(a_{i_{1}}\cdots a_{i_{j}})\cdot Z_{G},y=(a_{l_{1}}\cdots a_{l_{k}})\cdot Z_{G}$, then $yx=(a_{l_{1}}\cdots a_{l_{k}})(a_{i_{1}}\cdots a_{i_{j}})\cdot Z_{G}=(a_{i_{1}}\cdots a_{i_{j}})(a_{l_{1}}\cdots a_{l_{k}})d\cdot Z_{G}=(a_{i_{1}}\cdots a_{i_{j}})(a_{l_{1}}\cdots a_{l_{k}})\cdot Z_{G}=xy$ by Proposition \ref {5.3}.
\hfill $\Box$

\begin {Remark} \label {5.7} The converse of the above theorem is not true, and we give an example below.

Assume that $G=\{\scriptsize{\left (\begin{array} {lll}
1&x&y\\
0&1&z\\
0&0&1\\
\end {array} \right)}\mid x,y,z\in \mathbb R\}$ is a Heisenberg group, then $Z_{G}=\{\scriptsize{\left (\begin{array} {lll}
1&0&y\\
0&1&0\\
0&0&1\\
\end {array} \right)}$
$\mid y\in \mathbb R\}$. Consider the map $\phi:G\mapsto (\mathbb R\times \mathbb R,+)$ is given by $\scriptsize{\left (\begin{array} {lll}
1&x&y\\
0&1&z\\
0&0&1\\
\end {array} \right)}\mapsto (x,z)$, then $\phi$ is a group homomorphism from $G$ to abelian group $(\mathbb R\times \mathbb R,+)$ and $ker(\phi)=Z_{G}$. Thus $G/Z_{G}\cong (\mathbb R\times \mathbb R,+)$ is an abelian group, but $G$ is a non square commutative group since
$\scriptsize\left (\begin{array} {lll}
1&0&0\\
0&1&1\\
0&0&1\\
\end {array} \right)\left (\begin{array} {lll}
1&1&0\\
0&1&0\\
0&0&1\\
\end {array} \right))^{2}=\left (\begin{array} {lll}
1&2&1\\
0&1&2\\
0&0&1\\
\end {array} \right)\neq\left (\begin{array} {lll}
1&2&3\\
0&1&2\\
0&0&1\\
\end {array} \right)=(\left (\begin{array} {lll}
1&1&0\\
0&1&0\\
0&0&1\\
\end {array} \right)\left (\begin{array} {lll}
1&0&0\\
0&1&1\\
0&0&1\\
\end {array} \right))^{2}$.
\end {Remark}

\begin {Proposition} \label {5.8} Assume that $G$ is the group generated by $a,b$. Then $G$ is a square commutative group if and only if $G/Z_{G}$ is an abelian group and $(ab)^{2}=(ba)^{2}$.
\end {Proposition}
\noindent {\it Proof.} The only if part of the assertion is obvious by Proposition \ref {5.6}. It remains to prove the if part. Since $G/Z_{G}$ is an abelian group, then $ab=bad$, where $d\in Z_{G}$,
$(ab)^{2}=badbad=(ba)^{2}d^{2}$, so $d^{2}=e$, we have $a^{2}b=abad=ba^{2}d^{2}=ba^{2}$, $ab^{2}=badb=b^{2}ad^{2}=b^{2}a$, then $G$ is a square commutative group by Proposition \ref {1.3}.
\hfill $\Box$

The proof of Proposition \ref {5.8} shows that, if $ab=bad$, the other three can be deduced from any one of the four conditions $d^{2}=e, (ab)^{2}=(ba)^{2}, a^{2}b=ba^{2}, ab^{2}=b^{2}a$. Therefore, the proposition still holds if $(ab)^{2}=(ba)^{2}$ is replaced by one of $a^{2}b=ba^{2}$ and $ab^{2}=b^{2}a$.

For $x, y\in G$, if there exists $d\in Z_{G}^{2}$ such that  $x=yd$, then we write $x\sim y$. We know that $Z_{G}^{2}$ is a group and $\sim$ is an equivalent relation on $G$.

\begin {Theorem} \label {5.9} $G$ is a square commutative group if and only if $\widehat{G}$ is an abelian group. Meanwhile, we have $xy\sim yx$ for any $x,y\in G$.
\end {Theorem}
\noindent {\it Proof.} For any $x,y\in \widehat{G}$, we have $x=(a_{i_{1}}\cdots a_{i_{j}})\cdot Z_{G}^{2},y=(a_{l_{1}}\cdots a_{l_{k}})\cdot Z_{G}^{2}$, then $yx=(a_{l_{1}}\cdots a_{l_{k}})(a_{i_{1}}\cdots a_{i_{j}})\cdot Z_{G}^{2}=(a_{i_{1}}\cdots a_{i_{j}})(a_{l_{1}}\cdots a_{l_{k}})\cdot Z_{G}^{2}=xy$ by Proposition \ref {5.3}. Since $\widehat{G}$ is an abelian group, then $xy=yxd_{xy}$ and $d_{xy}^{2}=e$, $\forall\ x,y\in G,d_{xy}\in Z_{G}^{2}$, we have $(x_{1}x_{2})^{2}=x_{2}x_{1}d_{x_{1}x_{2}}x_{2}x_{1}d_{x_{1}x_{2}}
=(x_{2}x_{1})^{2}d_{x_{1}x_{2}}^{2}
=(x_{2}x_{1})^{2}$, $x_{1}x_{2}^{2}=x_{2}x_{1}d_{x_{1}x_{2}}x_{2}
=x_{2}^{2}x_{1}d_{x_{1}x_{2}}^{2}=x_{2}^{2}x_{1}$, $x_{1}(x_{2}x_{3})^{2}=x_{1}x_{2}x_{3}x_{2}x_{3}
=x_{2}x_{3}x_{2}x_{3}x_{1}d_{x_{1}x_{2}}^{2}d_{x_{3}x_{1}}^{2}
=(x_{2}x_{3})^{2}x_{1}$, where any $x_{1},x_{2},x_{3}$ are mutually exclusive elements in $T_{n}$, then $G$ is a square commutative group by Proposition \ref {1.3} and Proposition \ref {5.4}.
\hfill $\Box$

\begin {Corollary} \label {5.9'} Assume that $G$ is a p-group with $|G|=p^{n},p>2,n\in\mathbb N$. Then $G$ is a square commutative group if and only if $G$ is an abelian group.
\end {Corollary}
\noindent {\it Proof.} we obtain $Z_{G}^{2}=\{e\}$ since $d^{2}\neq e$ if $d\neq e$, then $\widehat{G}\cong G$.
\hfill $\Box$

\begin {Example} \label {5.10} Assume that $G$ is a p-group with $|G|=p^{3},p>2$. If $|Z_{G}|=p^{2}$, then $G/Z_{G}\cong \mathbb Z_{p}$ is a cyclic group, so $G$ is an abelian group, it is a contradiction. Consequently, $|Z_{G}|=p$ or $p^{3}$ since the center of the p-group must be non-trivial. If $|Z_{G}|=p^{3}$, then $G$ is an abelian group. If $|Z_{G}|=p$, then $|G/Z_{G}|=p^{2}$, $G/Z_{G}$ is an abelian group, but $G$ is an non-abelian group, thus $G$ is an non square commutative group.
\end {Example}

\begin {Example} \label {5.11} Let us consider $D_{2n}=\langle a,b\mid a^{n}=b^{2}=1,aba=b\rangle$ with $n=2k,k\in\mathbb N$. We have $Z_{D_{2n}}=\{e,a^{k}\}$ and $Z_{D_{2n}}^{2}=Z_{D_{2n}}$, then $\widehat{D_{2n}}=D_{2n}/Z_{D_{2n}}\cong D_{n}=D_{2k}$. It is known that $D_{2k}$ is a non-abelian group except for the case $k=2$. If $D_{2n}$ is a square commutative group, then $n=4$ by Theorem \ref {5.9}.
\end {Example}

\begin {Theorem} \label {5.12} Assume that $G$ is a group. Then $G^{[2]}\subset Z_{G}$ if and only if $G$ is a square commutative group.
\end {Theorem}
\noindent {\it Proof.} The if part of the assertion is obvious by Proposition \ref {5.3} and Theorem \ref {5.5}. It remains to prove the only if part. For $\forall\ x,y\in G$, we have $y^{2},(yx)^{2}\in Z_{G}$ since $G^{[2]}\subset Z_{G}$, then $xyxy^{2}=y^{2}xyx,(xy)^{2}y=y(yx)^{2}=(yx)^{2}y,(xy)^{2}=(yx)^{2}$.
\hfill $\Box$
\section{The square commutativity of $BS(p,q)$ with some
additional relations}
$BS(p,q)=\langle a,b\mid a^{p}b=ba^{q}\rangle$ is called the Baumslag-Solitar group (see \cite{L}), $p,q\in\mathbb Z^{*}$. Without any specific statement, the following $G$ denotes a group generated by $a,b$ with unit element $e$.

\begin {Lemma} \label {2.0} If $a$ is a generator satisfy the relation $a^{p}b=ba^{q}$ in $BS(p,q)$, then $ord(a)=\infty$.
\end {Lemma}
\noindent {\it Proof.} Let $G$ to be a matrix group generated by $\left (\scriptsize\begin{array} {lll}
1&1\\
0&1\\
\end {array} \right),\scriptsize\left (\begin{array} {lll}
\frac{p}{q}&0\\
0&1\\
\end {array} \right)$. Consider the map $\phi:BS(p,q)\mapsto G$ is given by $a\mapsto \scriptsize\left (\begin{array} {lll}
1&1\\
0&1\\
\end {array} \right),b\mapsto\left (\begin{array} {lll}
\frac{p}{q}&0\\
0&1\\
\end {array} \right)$, then $\phi$ is a group homomorphism from $BS(p,q)$ to $G$. If $ord(a)<\infty$, then $ord(\phi(a))<\infty$,
but $(\phi(a))^{n}=\scriptsize\left (\begin{array} {lll}
1&n\\
0&1\\
\end {array} \right)\neq\left (\begin{array} {lll}
1&0\\
0&1\\
\end {array} \right)$ for any $n\in \mathbb N$, which is a controdictory.
\hfill $\Box$

\begin {Proposition} \label {2.1}
{\rm (i)}  Assume that $G$ is a square commutative group and $a^{p}b=ba^{q}$. Then
$e= \left \{
\begin {array} {ll}
a^{2(p-q)},& \hbox {if} \hbox { $p,q$ } \hbox {is odd}\\
a^{p-q},& \hbox {if} \hbox { $p$ or $q$ } \hbox {is even}\\
\end {array} \right.$.

{\rm (ii)} Assuming that $a^{p}b=ba^{p}$ and $p$ is odd, then $G$ is a square commutative group if and only if $G$ is an abelian group.
\end {Proposition}
\noindent {\it Proof.} {\rm (i)} The conclusion is clearly valid when $p=q$. Assume that $p\neq q$. Since $G$ is a square commutative group, $a^{2}b=ba^{2}$ by Lemma \ref {1.2}. If $p$ is even, $a^{p}b=ba^{p}=ba^{q},a^{p-q}=e$; if $q$ is even, $a^{p}b=a^{q}b=ba^{q},a^{p-q}=e$; if $p,q$ are both odd, $a^{p}b=aba^{p-1}=ba^{q},ab=ba^{q-p+1}$, then $a^{2}b=ba^{2q-2p+2}=ba^{2}$, which implys $a^{2(p-q)}=e$.

{\rm (ii)} If $p$ is odd and $G$ is a square commutative group, then $a^{p}b=aba^{p-1}=ba^{p},ab=ba$, and $G$ is an abelian group by Proposition \ref {1.1}{\rm (ii)}.
\hfill $\Box$

When $p\neq q$, if $BS(p,q)$ is the square commutative group, then $a$ is finite order by Proposition \ref {2.1}{\rm (i)}, which is contradictory to Lemma \ref {2.0}, so only $p=q$ is to be considered. If $p=1$, then $BS(p,p)$ is an abelian group. If $p>2$, by Lemma \ref {1.2}, $ba^{2}=a^{2}b$, hence $p=2$, which contradicts with $p>2$, so $BS(p,p)$ is not a square commutative group. If $p=2$, then $BS(2,2)$ is not a square commutative group since $ba^{2}=a^{2}b$ can not imply $ab^{2}=b^{2}a$, if $BS(2,2)$ also satisfy the relation $ab^{2}=b^{2}a$, we add additional conditions, this relation should still hold. The following is the counter example, $G=\langle a,b\mid a^{4}=b^{3}=e,ba=ab^{2},ba^{2}=a^{2}b\rangle=
\{e,a,a^{2},a^{3},b,ab,a^{2}b,a^{3}b,b^{2},ab^{2},a^{2}b^{2},a^{3}b^{2}
\}=
\{e,a,a^{2},a^{3}\}\cdot\{e,b,b^{2}\}$.

Synthesizing the above argument, we have the following proposition.

\begin {Proposition} \label {2.1'}
$BS(p,q)$ is a square commutative group if and only if $p=q=1$.
\end {Proposition}

The previous example was given where $ba^{2}=a^{2}b,(ab)^{2}=(ba)^{2}$ is valid but $ab^{2}=b^{2}a$ is not, and below we will give an example where $ba^{2}=a^{2}b,ab^{2}=b^{2}a$ is valid while $(ab)^{2}=(ba)^{2}$ is not, e.g., $G=\langle a,b\mid a^{4}=b^{2}=e,ba=(ab)^{3},ba^{2}=a^{2}b\rangle=
\{e,a,a^{2},a^{3},b,ab,a^{2}b,a^{3}b,$
$ba,aba,a^{2}ba,
a^{3}ba,bab,abab,a^{2}bab,a^{3}bab\}=
\{e,a,a^{2},a^{3}\}\cdot\{e,ba\}\cdot\{e,b\}$.

\begin {Theorem} \label {2.2} Assuming $ab=ba^{j}$, $j\in \mathbb Z^{*}$, the necessary and sufficient condition for $G$ to be a square commutative group is
$e= \left \{
\begin {array} {ll}
a^{2(j-1)}& \hbox {if} \hbox { $j$ } \hbox {is odd}\\
a^{j-1}& \hbox {if} \hbox { $j$ } \hbox {is even}
\end {array} \right.$.
\end {Theorem}

\noindent {\it Proof.} By Proposition \ref {2.1} , it is sufficient to prove sufficiency. When $j$ is even, $ab=ba^{j}=ba$, and the conclusion is obvious. Now assume that $j$ is odd, and that $a^{2}b=ba^{2j}=ba^{2}$. Using induction on $k\in \mathbb N_{0}$, we get
$a^{k}b= \left \{
\begin {array} {ll}
ba^{(j-1)+k}& \hbox {if} \hbox { $k$ } \hbox {is odd}\\
ba^{k}& \hbox {if} \hbox { $k$ } \hbox {is even}
\end {array} \right.$. Use induction on $k,l\in \mathbb N_{0}$, we get
$a^{k}b^{l}= \left \{
\begin {array} {ll}
b^{l}a^{(j-1)+k}& \hbox {if} \hbox { $k,l$ } \hbox {is odd}\\
b^{l}a^{k}& \hbox {if} \hbox { $k$ or $l$ } \hbox {is even}
\end {array} \right.$.  Take $\forall\ b^{l_{1}}a^{k_{1}},b^{l_{2}}a^{k_{2}}\in G$, then

$b^{l_{1}}a^{k_{1}}b^{l_{2}}a^{k_{2}}= \left \{
\begin {array} {ll}
b^{l_{1}+l_{2}}a^{k_{1}+k_{2}+(j-1)}& \hbox {if} \hbox { $k_{1},l_{2}$ } \hbox {is odd}\\
b^{l_{1}+l_{2}}a^{k_{1}+k_{2}}& \hbox {if} \hbox { $k_{1}$ or $l_{2}$ } \hbox {is even}
\end {array} \right.$.

$b^{l_{2}}a^{k_{2}}b^{l_{1}}a^{k_{1}}= \left \{
\begin {array} {ll}
b^{l_{1}+l_{2}}a^{k_{1}+k_{2}+(j-1)}& \hbox {if} \hbox { $k_{2},l_{1}$ } \hbox {is odd}\\
b^{l_{1}+l_{2}}a^{k_{1}+k_{2}}& \hbox {if} \hbox { $k_{2}$ or $l_{1}$ } \hbox {is even}
\end {array} \right.$.

$b^{l}a^{k}b^{l}a^{k}= \left \{
\begin {array} {ll}
b^{2l}a^{(j-1)+2k}& \hbox {if} \hbox { $k,l$ } \hbox {is odd}\\
b^{2l}a^{2k}& \hbox {if} \hbox { $k$ or $l$ } \hbox {is even}
\end {array} \right.$.

Then $(b^{l_{1}}a^{k_{1}}b^{l_{2}}a^{k_{2}})^{2}= \left \{
\begin {array} {ll}
b^{l_{1}+l_{2}}a^{k_{1}+k_{2}+(j-1)}b^{l_{1}+l_{2}}a^{k_{1}+k_{2}+(j-1)}& \hbox {if} \hbox { $k_{1},l_{2}$ } \hbox {is odd}\\
b^{l_{1}+l_{2}}a^{k_{1}+k_{2}}b^{l_{1}+l_{2}}a^{k_{1}+k_{2}}& \hbox {if} \hbox { $k_{1}$ or $l_{2}$ } \hbox {is even}\\
\end {array} \right. \\
=\left \{
\begin {array} {ll}
b^{2(l_{1}+l_{2})}a^{2(k_{1}+k_{2})+(j-1)}& \hbox {if} \hbox { $k_{1},l_{2}$ } \hbox {is odd,} \hbox { $k_{2},l_{1}$ } \hbox {is even}\\
b^{2(l_{1}+l_{2})}a^{2(k_{1}+k_{2})}& \hbox {if} \hbox { $k_{1},l_{2}$ } \hbox {is odd,}\hbox { $k_{2}$ or $l_{1}$ } \hbox {is odd}\\
b^{2(l_{1}+l_{2})}a^{2(k_{1}+k_{2})+(j-1)}& \hbox {if} \hbox { $k_{1}$ or $l_{2}$ } \hbox {is even,}\hbox { $k_{1}+k_{2},l_{1}+l_{2}$ } \hbox {is odd}\\
b^{2(l_{1}+l_{2})}a^{2(k_{1}+k_{2})}& \hbox {if} \hbox { $k_{1}$ or $l_{2}$ } \hbox {is even,}\hbox { $k_{1}+k_{2}$ or $l_{1}+l_{2}$ } \hbox {is even}\\
\end {array} \right.\\
=\left \{
\begin {array} {ll}
b^{2(l_{1}+l_{2})}a^{2(k_{1}+k_{2})+(j-1)}& \hbox {if} \hbox { $k_{1}+k_{2},l_{1}+l_{2}$ } \hbox {is odd}\\
b^{2(l_{1}+l_{2})}a^{2(k_{1}+k_{2})}& \hbox {if}\hbox { $k_{1}+k_{2}$ or $l_{1}+l_{2}$ } \hbox {is even}\\
\end {array} \right.$.

$(b^{l_{2}}a^{k_{2}}b^{l_{1}}a^{k_{1}})^{2}= \left \{
\begin {array} {ll}
b^{l_{1}+l_{2}}a^{k_{1}+k_{2}+(j-1)}b^{l_{1}+l_{2}}a^{k_{1}+k_{2}+(j-1)}& \hbox {if} \hbox { $k_{2},l_{1}$ } \hbox {is odd}\\
b^{l_{1}+l_{2}}a^{k_{1}+k_{2}}b^{l_{1}+l_{2}}a^{k_{1}+k_{2}}& \hbox {if} \hbox { $k_{2}$ or $l_{1}$ } \hbox {is even}\\
\end {array} \right. \\
=\left \{
\begin {array} {ll}
b^{2(l_{1}+l_{2})}a^{2(k_{1}+k_{2})+(j-1)}& \hbox {if} \hbox { $k_{2},l_{1}$ } \hbox {is odd,} \hbox { $k_{1},l_{2}$ } \hbox {is even}\\
b^{2(l_{1}+l_{2})}a^{2(k_{1}+k_{2})}& \hbox {if} \hbox { $k_{2},l_{1}$ } \hbox {is odd,}\hbox { $k_{1}$ or $l_{2}$ } \hbox {is odd}\\
b^{2(l_{1}+l_{2})}a^{2(k_{1}+k_{2})+(j-1)}& \hbox {if} \hbox { $k_{2}$ or $l_{1}$ } \hbox {is even,}\hbox { $k_{1}+k_{2},l_{1}+l_{2}$ } \hbox {is odd}\\
b^{2(l_{1}+l_{2})}a^{2(k_{1}+k_{2})}& \hbox {if} \hbox { $k_{2}$ or $l_{1}$ } \hbox {is even,}\hbox { $k_{1}+k_{2}$ or $l_{1}+l_{2}$ } \hbox {is even}\\
\end {array} \right.\\
=\left \{
\begin {array} {ll}
b^{2(l_{1}+l_{2})}a^{2(k_{1}+k_{2})+(j-1)}& \hbox {if} \hbox { $k_{1}+k_{2},l_{1}+l_{2}$ } \hbox {is odd}\\
b^{2(l_{1}+l_{2})}a^{2(k_{1}+k_{2})}& \hbox {if}\hbox { $k_{1}+k_{2}$ or $l_{1}+l_{2}$ } \hbox {is even}\\
\end {array} \right.$.

Then $(b^{l_{1}}a^{k_{1}}b^{l_{2}}a^{k_{2}})^{2}=
(b^{l_{2}}a^{k_{2}}b^{l_{1}}a^{k_{1}})^{2}$ holds for all $k_{1},k_{2},l_{1},l_{2}\in \mathbb N_{0}$, hence $G$ is the square commutative group. \hfill $\Box$

\begin {Corollary} \label {1.9} {\rm (i)} Suppose $G=\langle a,b\mid a^{n}=e,ab=ba^{j}\rangle$, and $2\leq n,j\in \mathbb N,j<n$, then $G$ is a square commutative group if and only if $j$ is odd and $n=2(j-1)$.

{\rm (ii)} Assuming that $G=\langle a,b\mid a^{n}=b^{m}=e,ab=ba^{j}\rangle$,
where $2\leq n,m,j\in \mathbb N,j<n$, then $G=\{b^{i}a^{j}\mid 0\leq i\leq m-1,0\leq j\leq n-1\}$ and $\mid G\mid =mn$. Assume that $G$ is a square commutative group, so that $n=2(j-1)$, and if $j=n-1$, then $n=4$. Further, if $m=2$, then $G$ is a dihedral group, $\mid G\mid =2n$, which is an another proof of \cite[Proposition 1]{W}. At this point it can be seen that $\langle a,b\mid a^{4}=b^{2}=e,aba=b\rangle\cong\langle a,b\mid a^{2}=b^{2}=e,(ab)^{2}=(ba)^{2}\rangle$.

{\rm (iii)} Assuming that $a,b$ satisfy the relation $a^{p}b=ba^{q},b^{2}=e,|p|,|q|\geq2$, then $G$ is the square commutative group if and only if $a,b$ satisfy the relation $(ab)^{2}=(ba)^{2},p=q=\pm2$.
\end {Corollary}

\begin {Proposition} \label {2.10} Assume that $a^{2}b=ba^{2}$. Then $G$ is a square commutative group if and only if $G/Z_{G}$ is an abelian group.
\end {Proposition}
\noindent {\it Proof.} It is clear by Proposition \ref {5.8}.
\hfill $\Box$

Weicai Wu\\
School of Mathematics and Statistics, Guangxi Normal University,
Guilin 541004, Guangxi, People's Republic of
China.\\
E-mail: \textsf{weicaiwu@hnu.edu.cn}\\[0.3cm]
Mingxuan Yang\\
School of Mathematics and Statistics, Guangxi Normal University,
Guilin 541004, Guangxi, People's Republic of
China.\\
E-mail: \textsf{3223237115@qq.com}\\[0.3cm]
Yangbo Zhou\\
Changsha Xingsha Middle School,Xingsha 410134,  Changsha, People's Republic of
China.\\
E-mail: \textsf{yangbozhou2024@163.com}\\[0.3cm]
Chao Rong\\
School of Mathematics and Statistics, Guangxi Normal University,
Guilin 541004, Guangxi, People's Republic of
China.\\
E-mail: \textsf{2207318475@qq.com}
\end{document}